\documentclass[a4paper,3p,11pt,sort&compress,review]{elsarticle}

\RequirePackage{snapshot}

\usepackage[Symbol]{upgreek}
\usepackage[framed,numbered,autolinebreaks,useliterate]{mcode}
\usepackage[T1]{fontenc}
\usepackage[latin9]{inputenc}
\usepackage{amsmath}
\usepackage{amssymb}
\usepackage{graphicx}
\usepackage{setspace}
\usepackage{esint}
\PassOptionsToPackage{normalem}{ulem}
\usepackage{ulem}
\usepackage{marvosym}
\usepackage{multirow}
\usepackage{makecell}
\usepackage{soul}
\usepackage{physics}
\usepackage{xfrac}
\usepackage{dsfont}
\usepackage[hidelinks]{hyperref}
\usepackage{nomencl}
\usepackage{adjustbox}
\usepackage{mathtools}
\usepackage{nicefrac}
\usepackage[hang,flushmargin]{footmisc}
\usepackage{amsthm}
\usepackage{lineno}
\usepackage{units}

\usepackage{enumitem}
\usepackage{subfig}
\usepackage{bm}
\usepackage{color}
\usepackage{epstopdf}

\DeclareMathAlphabet\mathbfcal{OMS}{cmsy}{b}{n}      
\DeclareMathAlphabet\mathcal{OMS}{cmsy}{m}{n}      

\setlist{nolistsep}

\DeclareMathAlphabet\mathbfcal{OMS}{cmsy}{b}{n}

\makeatletter
\newcommand{\shorteq}{%
  \settowidth{\@tempdima}{-}
  \resizebox{\@tempdima}{\height}{=}%
}
\makeatother

\makeatletter
\renewcommand*\env@matrix[1][*\c@MaxMatrixCols c]{%
\hskip -\arraycolsep
\let\@ifnextchar\new@ifnextchar
\array{#1}}
\makeatother

\DeclareMathSymbol{\shortminus}{\mathbin}{AMSa}{"39}

\soulregister\cite7
\soulregister\ref7
\soulregister\pageref7

\usepackage{hyperref}

\makeatletter

\newcommand{\Her}[0]{{\mathsf{H}}}                  
\newcommand{\I}[0]{\mathrm{i}\mkern1mu}             
\newcommand{\Ma}[0]{\mathbfcal{M}}                  
\newcommand{\La}[0]{\mathbfcal{L}}                  

\newcommand*{\noi}{\noindent}

\DeclareMathAlphabet{\mathup}{OT1}{\familydefault}{m}{n}

\newcommand*\xbar[1]{%
  \hbox{%
    \vbox{%
      \hrule height 0.5pt 
      \kern0.3ex
      \hbox{%
        \kern-0.1em
        \ensuremath{#1}%
        \kern-0.1em
      }%
    }%
  }%
}

\newsavebox\myboxA
\newsavebox\myboxB
\newlength\mylenA

\newcommand*\xoverline[2][0.75]{%
    \sbox{\myboxA}{$\m@th#2$}%
    \setbox\myboxB\null
    \ht\myboxB=\ht\myboxA%
    \dp\myboxB=\dp\myboxA%
    \wd\myboxB=#1\wd\myboxA
    \sbox\myboxB{$\m@th\overline{\copy\myboxB}$}
    \setlength\mylenA{\the\wd\myboxA}
    \addtolength\mylenA{-\the\wd\myboxB}%
    \ifdim\wd\myboxB<\wd\myboxA%
       \rlap{\hskip 0.5\mylenA\usebox\myboxB}{\usebox\myboxA}%
    \else
        \hskip -0.5\mylenA\rlap{\usebox\myboxA}{\hskip 0.5\mylenA\usebox\myboxB}%
    \fi}

\usepackage{tikz,pgfplots,grffile} 
\pgfplotsset{compat=newest}
\usetikzlibrary{calc, shapes, arrows}
\tikzstyle{arrow} = [thick,->,>=stealth]

\graphicspath{{figures/}}

\definecolor{myblue}{rgb}{0.3804, 0.6471, 0.7608}

\journal{}

\newdefinition{rmk}{Remark}

\makeatother
\PassOptionsToPackage{numbers,sort&compress}{natbib}

\makenomenclature
\begin{document}
\makeatletter
\def\bm@pmb@#1{{%
      \setbox\tw@\hbox{$\m@th\mkern1.5mu$}%
      \mathchoice
      \bm@pmb@@\displaystyle\@empty{#1}%
      \bm@pmb@@\textstyle\@empty{#1}%
      \bm@pmb@@\scriptstyle\defaultscriptratio{#1}%
      \bm@pmb@@\scriptscriptstyle\defaultscriptscriptratio{#1}}}
\makeatother

\title{Computing critical velocities in waveguides via multiparameter eigenvalue problems}

\author[ovgu]{Hauke~Gravenkamp\corref{cor1}}
\ead{hauke.gravenkamp@ovgu.de}

\author[csu]{Xiang~Liu}

\address[ovgu]{Institute of Materials, Technologies and Mechanics, Otto von Guericke University Magdeburg\\ 39106 Magdeburg, Germany}

\address[csu]{School of Traffic and \& Transportation Engineering, Central South University\\ Changsha, China}

\cortext[cor1]{Corresponding author}

\begin{abstract}
  In waveguide dynamics and moving-load problems (e.g., high-speed trains), critical velocities indicate the onset of strong vibration amplification. In systems that are invariant in the direction of motion, these velocities can be identified from dispersion relations as points where the phase and group velocities of a propagating mode coincide. Finding such points indirectly by tracing dispersion curves can be cumbersome and potentially unreliable for multimodal systems with complex branch interactions.
  We present a direct method for computing critical velocities in such scenarios, specifically in the context of semi-analytical methods. Starting from a polynomial parameter-dependent eigenvalue problem for the wavenumber--frequency relation, incorporating the additional condition of equal phase and group velocities yields a singular polynomial multiparameter eigenvalue problem that can be linearized and solved using established algorithms.
  The proposed approach enables the simultaneous computation of all critical points without requiring the tracing of dispersion curves. Its performance is demonstrated by several benchmark problems, confirming the accurate and robust identification of critical velocities.
\end{abstract}
\begin{keyword}
  critical velocity; wave propagation;  multiparameter eigenvalue problem; dispersion; high-speed trains
\end{keyword}
\maketitle


\section{Introduction}\label{sec:intro}\noi
Critical velocities are a central concept in high-speed transportation and, more generally, in moving-load and waveguide dynamics \cite{Bierer2007, Kausel2020}. In railway engineering, the phenomenon is commonly associated with a rapid amplification of track and ground vibration when the train speed approaches a characteristic wave speed of the coupled track--embankment--soil system. Early theoretical work already linked this amplification to the generation of strong surface-wave radiation by superfast trains~\cite{Krylov1995}. More recent studies have established that the relevant threshold is governed by the dispersive wave-propagation characteristics of the supporting system and, in particular, by minima of the phase-velocity spectrum rather than by a classical resonance of a finite structure \cite{Mezher2016, Kausel2020, Estaire2024}. Specifically, critical velocities occur when the phase velocity $c_p$ of a propagating mode equals its group velocity $c_g$.
The same physical idea reappears in overhead contact systems, where the so-called catenary barrier reflects the interaction between the contact-point speed and the wave-propagation properties of the tensioned cable system \cite{Vesali2021,Liu2026}. It also extends naturally to high-speed magnetically levitated vehicles, where wave-induced instability becomes relevant once the operating speed enters the supercritical regime \cite{Faragau2025}. Beyond transportation, the broader mechanics community has studied how wave speeds can be tailored or even self-controlled in nonlinear media, which further underlines the relevance of robust tools for locating characteristic propagation thresholds \cite{Xin2017}.

From a modeling perspective, many dynamic systems of interest are invariant, or at least locally periodic, in the direction of motion. Such structures include beams on elastic or viscoelastic foundations, layered soil profiles, periodically supported rails, overhead contact lines, and general prismatic structures. Such systems are naturally described in terms of dispersion relations between frequency and wavenumber. In the railway context, this viewpoint has been employed for critical-speed prediction in coupled track--ground models \cite{Mezher2016, Kausel2020, Estaire2024} and was recently extended to periodically varying viscoelastic foundations, where Floquet-type arguments become essential for determining the critical train speed \cite{Tran2025}. In parallel, semi-analytical and dynamic-stiffness-based waveguide models have matured considerably. These include approaches based on the thin layer method (TLM) \cite{Kausel1981a, Kausel1994}, the scaled boundary finite element method (SBFEM) \cite{Gravenkamp2012e, Gravenkamp2018}, and semi-analytical finite elements (SAFE) \cite{Marzani2008, Hayashi2003a}. These three methods all employ a numerical discretization of a waveguide's cross-section, typically based on conventional finite elements, spectral elements, or splines \cite{Gravenkamp2021b}.

In contrast, the dynamic-stiffness method employs exact dispersion relations in the wavenumber-frequency domain for periodic structures and plate-built waveguides. Here, it was shown that the Wittrick--Williams algorithm can be adapted and enhanced for exact vibration and wave-propagation analyses, including explicit treatment of the $J_0$ count, which is crucial for robust and efficient computations over broad frequency ranges \cite{Liu2022a, Liu2024b, Zhou2024, Liu2026}. A common denominator of all the aforementioned approaches is that the underlying physics can be reduced, after semi-discretization in the cross-section, to a parameter-dependent eigenvalue problem of the form $\mathbf{W}(k,\omega)\,\mathbf{u} = \mathbf{0} $, where the wavenumber $k$ is the eigenvalue and the frequency $\omega$ plays the role of an independent parameter, or vice versa. It is worth noting that the dynamic stiffness method generally leads to very small \textit{transcendental} eigenvalue problems, whereas the semi-analytical methods yield somewhat larger but \textit{polynomial} eigenvalue problems. Hence, both classes of methods require very different solution procedures and entail distinct benefits and drawbacks. In this paper, we focus on semi-analytical methods and exploit their simpler polynomial structure.

From an algorithmic point of view, critical velocities are usually still obtained indirectly. The standard procedure relies on computing dispersion curves by solving for $k$ at a sequence of frequencies (or vice versa) to trace the branches, and then determining the critical points by searching for locations that satisfy the condition $c_p=c_g$. For multimodal systems, this approach is cumbersome and potentially unreliable, because dispersion diagrams may exhibit crossings, osculations, and veering. These issues are well-documented in semi-analytical waveguide models \cite{Kausel2015, Gravenkamp2023i}. A closely related problem is the computation of zero-group-velocity (ZGV) points in waveguides, i.e., points on the dispersion curves with $c_g=0$. It has recently been shown that such ZGV points can be computed directly as the solution of a multiparameter eigenvalue problem \cite{Kiefer2022}.
In a nutshell, the parameter-dependent eigenvalue problem mentioned earlier, together with the condition imposed on the group velocity, yields two coupled polynomial eigenvalue problems for the wavenumbers and frequencies at ZGV points. Such multiparameter eigenvalue problems are well understood \cite{Muhic2010}, and numerical methods for the solution are available \cite{Plestenjak2023}. Other related applications of this framework were found in the direct computation of leaky waves in layered structures coupled to unbounded media \cite{Gravenkamp2025, Gravenkamp2024a} as well as in the computation of Hopf bifurcations in fluid-conveying pipes \cite{Gravenkamp2026}.

Against this background, we will show in this paper that the approach outlined in \cite{Kiefer2022} for ZGV points can be adapted to the computation of critical velocities, which corresponds to incorporating the condition $c_g=c_p$ (rather than $c_g = 0$) into a multiparameter eigenvalue problem. This method allows the direct computation of all critical points in a single step rather than as a posteriori information extracted from traced dispersion branches. The resulting formulation is sufficiently general to cover simple benchmark systems such as a beam on an elastic foundation, layered media representing track support, and more general three-dimensional waveguides, while at the same time opening the door to future extensions to periodic structures, leaky-wave settings, and large coupled systems.

\section{Problem statement}\label{sec:problem}
\begin{figure}\centering
  \includegraphics[width = 0.6\textwidth]{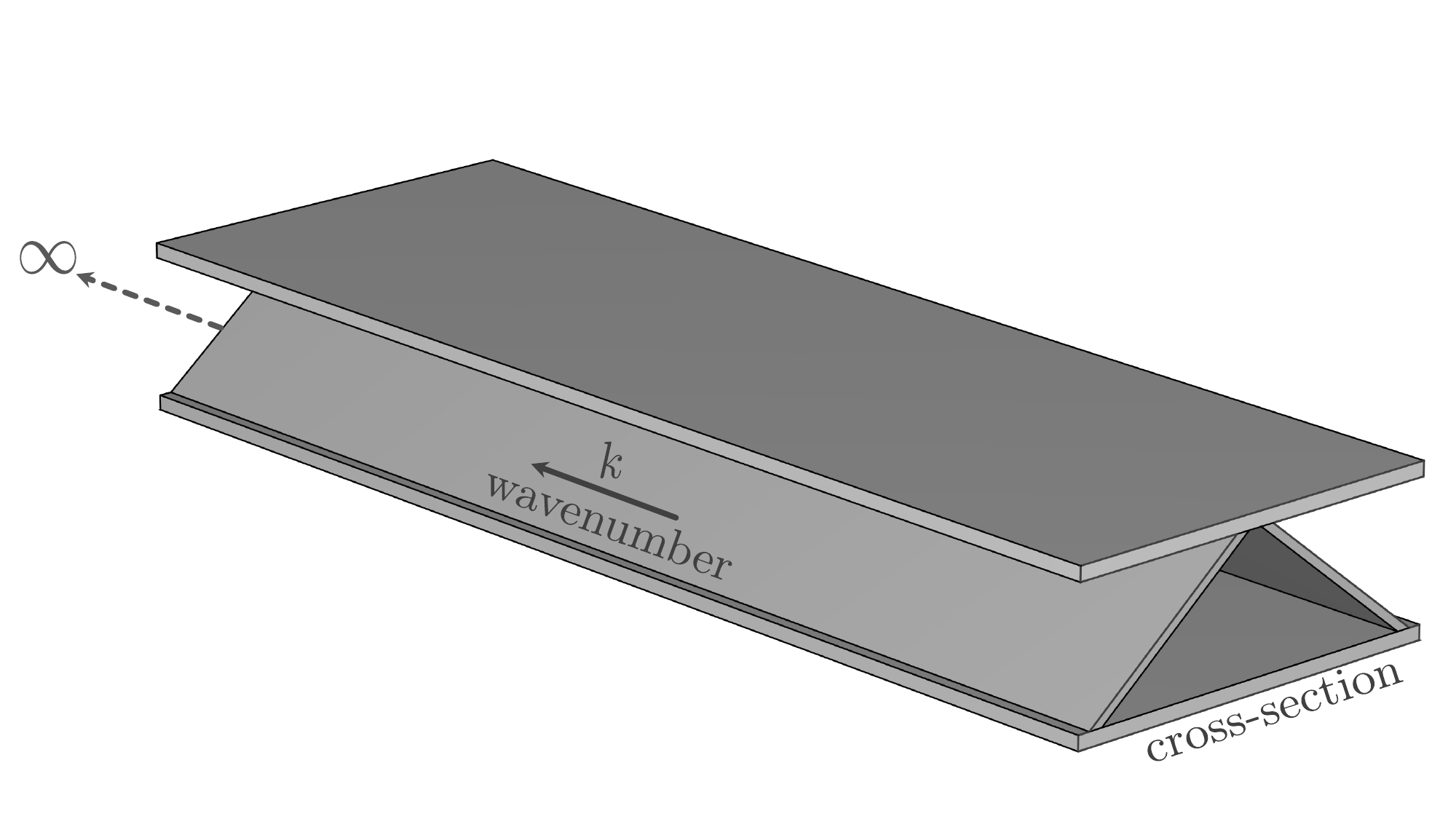}
  \caption{Example of a prismatic structure of uniform cross-section and infinite extent. \label{fig:prismaticGeometry}}
\end{figure}

\noi
We consider elastic waves propagating along infinitely long structures of uniform cross-section. A classic example is given by guided waves in layered plates of constant thickness, usually referred to as Lamb waves \cite{Lamb1917}. However, the cross-section can generally be of any shape \cite{Gravenkamp2013e}, see Fig.~\ref{fig:prismaticGeometry} for an example. Wave motion in such structures is characterized by modes, each propagating with a specific frequency-dependent axial wavenumber $k$ as well as a mode shape $\mathbf{u}$ that describes the displacement amplitudes on the cross-section.

Generally, $k$ and $\mathbf{u}$ are obtained as solutions to a parameter-dependent eigenvalue problem
\begin{equation}
  \label{eq:evp_discrete}
  \mathbf{W}(k, \omega)\, \mathbf{u} = \mathbf{0},
\end{equation}
where the frequency $\omega$ plays the role of an independent parameter. Solving this eigenvalue problem for the wavenumbers $k$ at varying $\omega$ results in continuous curves $k(\omega)$ for each mode, referred to as dispersion curves.
Here, $\mathbf{W}$ is a parameterized $n\negthinspace \times\negthinspace n$-matrix whose properties depend on the employed modeling approach, typically involving numerical approximations. A particularly relevant class of methods discretizes the cross-section by finite elements or closely related numerical methods, resulting in a quadratic matrix function of the form \cite{Kausel1981a, Hayashi2003a, Marzani2008, Gravenkamp2012e}
\begin{equation}\label{eq:operator_discrete}
  \mathbf{W}(k, \omega) = -k^2 \mathbf{L}_2 + \I k \mathbf{L}_1 + \mathbf{L}_0 + \omega^2 \mathbf{M}
\end{equation}
with real matrices $\mathbf{L}_i$, $\mathbf{M}$. We will also briefly discuss the following variant of a fourth-order eigenvalue problem that emerges in modeling, e.g., beams on an elastic foundation \cite{Kausel2020}:
\begin{equation}\label{eq:operator_discrete_order4}
  \mathbf{W}(k, \omega) = k^4 \mathbf{L}_4 + \mathbf{L}_0 + \omega^2 \mathbf{M}.
\end{equation}
Note that this problem is linear in $k^4$ and $\omega^2$. The formulations leading to Eqs.~\eqref{eq:operator_discrete} and \eqref{eq:operator_discrete_order4} are well established in the literature and do not require repetition (see the aforementioned references for details).
Here, we are interested in the critical velocities, which are isolated points on the dispersion curves where the phase velocity equals the group velocity, i.e.,
\begin{equation}\label{eq:cp_equals_cg}
  c_p \coloneq \frac{\omega}{k} = \frac{\partial \omega}{\partial k} \eqcolon c_g, \qquad \omega\in\mathbb{R},\quad k \in \mathbb{R}\negthinspace\setminus\negthinspace\{0\}.
\end{equation}
Note that, for the phase and group velocities to be well-defined, we assume real-valued frequencies and wavenumbers, which is true for propagating modes in lossless media.
For a given eigenvector/eigenvalue pair corresponding to Eq.~\eqref{eq:operator_discrete} and assuming a Hermitian structure\footnote{The most common semi-analytical formulations lead to a Hermitian eigenvalue problem. There are, however, relevant scenarios where this property does not hold, in particular, some models for material damping, non-symmetric formulations of acoustic/elastic coupling, or cross-section discretizations by Petrov-Galerkin methods.}, the group velocity can be computed a posteriori as
\begin{equation}\label{eq:groupVel}
  c_g = \frac{\mathbf{u}^\Her (2k\mathbf{L}_2 - \I \mathbf{L}_1)\mathbf{u}}{2\omega\,\mathbf{u}^\Her\mathbf{M}\mathbf{u}}.
\end{equation}

\section{Computation of critical velocities}\label{sec:cCrit}\noi
To derive a formulation for the direct computation of critical velocities, we will differentiate the eigenvalue problem \eqref{eq:evp_discrete} with respect to the wavenumber and substitute the condition \eqref{eq:cp_equals_cg}. The resulting equation, together with the original eigenvalue problem, poses two coupled eigenvalue problems (i.e., a two-parameter eigenvalue problem), whose solutions contain the sought critical points.
The basic formulation is similar to that previously presented for the computation of ZGV points \cite{Kiefer2022, Plestenjak2025}, except for a slightly different structure arising from an additional term due to the nonzero group velocity. Note that, in the following, the \textit{dash} symbol always refers to derivatives with respect to the wavenumber $k$.

We begin with the matrix operator defined by Eq.~\eqref{eq:operator_discrete}.
Differentiating the eigenvalue problem with respect to the wavenumber yields
\begin{equation}
  \mathbf{W}'\, \mathbf{u} + \mathbf{W}\, \mathbf{u}' = \mathbf{0}
\end{equation}
with
\begin{equation}
  \mathbf{W}'(k, \omega) = -2k\,\mathbf{L}_2 + \I \mathbf{L}_1 + 2\omega\,\omega'\,\mathbf{M}.
\end{equation}
Note that $\omega'$ is the group velocity.
At a critical point, we substitute the condition given by Eq.~\eqref{eq:cp_equals_cg} and define the resulting matrix function as
\begin{equation}
  \mathbf{W}_c(k, \omega) \coloneqq -2k\,\mathbf{L}_2 + \I \mathbf{L}_1 + 2\,\frac{\omega^2}{k}\,\mathbf{M}.
\end{equation}
Hence, a solution of the eigenvalue problem with the additional property $c_g=c_p$ must satisfy $\mathbf{W}_c\, \mathbf{u} + \mathbf{W}\, \mathbf{u}' = \mathbf{0}$ as well as $\mathbf{W}\, \mathbf{u}  = \mathbf{0}$, or in block matrix form:
\begin{equation}\label{eq:criticalPointSystem}
  \begin{bmatrix}
    \mathbf{W}   & \mathbf{0} \\
    \mathbf{W}_c & \mathbf{W}
  \end{bmatrix}
  \begin{bmatrix}
    \mathbf{u} \\ \mathbf{u}'
  \end{bmatrix} =
  \begin{bmatrix}
    \mathbf{0} \\ \mathbf{0}
  \end{bmatrix}.
\end{equation}
We multiply the second equation by $k$ to separate the parameters:
\begin{equation}\label{eq:criticalPointSystemK}
  \begin{bmatrix}
    \mathbf{W}    & \mathbf{0} \\
    k\mathbf{W}_c & \mathbf{W}
  \end{bmatrix}
  \begin{bmatrix}
    \mathbf{u} \\ k\mathbf{u}'
  \end{bmatrix}
  =
  \begin{bmatrix}
    \mathbf{0} \\ \mathbf{0}
  \end{bmatrix}.
\end{equation}
Substituting the explicit expressions for $\mathbf{W}$ and $\mathbf{W}_c$ yields
\begin{equation}\label{eq:criticalPointSystemExplicit}
  \left(
  -k^2
  \begin{bmatrix}
      \mathbf{L}_2  & \mathbf{0}   \\
      2\mathbf{L}_2 & \mathbf{L}_2
    \end{bmatrix}
  + \I k
  \begin{bmatrix}
      \mathbf{L}_1 & \mathbf{0}   \\
      \mathbf{L}_1 & \mathbf{L}_1
    \end{bmatrix}
  +
  \begin{bmatrix}
      \mathbf{L}_0 & \mathbf{0}   \\
      \mathbf{0}   & \mathbf{L}_0
    \end{bmatrix}
  + \omega^2
  \begin{bmatrix}
      \mathbf{M}  & \mathbf{0} \\
      2\mathbf{M} & \mathbf{M}
    \end{bmatrix}
  \right)
  \begin{bmatrix}
    \mathbf{u} \\ k\mathbf{u}'
  \end{bmatrix}
  =
  \begin{bmatrix}
    \mathbf{0} \\ \mathbf{0}
  \end{bmatrix},
\end{equation}
which we abbreviate as
\begin{equation}\label{eq:criticalPointSystemShort}
  \big(-k^2 \La_2 + \I k \La_1 + \La_0 + \omega^2 \Ma\big) \mathbf{v} = \mathbf{0}.
\end{equation}
The above equation constitutes an additional condition for an eigenvalue solution to be a critical point. This equation, together with the original eigenvalue problem \eqref{eq:evp_discrete} form a system of two eigenvalue problems in $\omega$ and $k$, whose solutions correspond to critical points. In other words, we obtain a \textit{quadratic} two-parameter eigenvalue problem:
\begin{subequations}
  \label{eq:2PEVP}
  \begin{align}
    (-k^2 \mathbf{L}_2 + \I k \mathbf{L}_1 + \mathbf{L}_0 + \omega^2 \mathbf{M})\mathbf{u} & = \mathbf{0},\label{eq:2ep1} \\
    (-k^2 \La_2 + \I k \La_1 + \La_0 + \omega^2 \Ma) \mathbf{v}                            & = \mathbf{0}. \label{eq:2ep2}
  \end{align}
\end{subequations}
For simplification and conciseness, it is rewritten as a \textit{linear} three-parameter eigenvalue problem, defining the parameters $\eta = (\I k)^2$, $\lambda = \I k$, $\mu = \omega^2$ and including a third equation to define the relationship between $\eta$ and $\lambda$:
\begin{subequations}
  \label{eq:3PEVP}
  \begin{align}
    (\eta \mathbf{L}_2 + \lambda \mathbf{L}_1 + \mathbf{L}_0 + \mu \mathbf{M})\mathbf{u} & = \mathbf{0},\label{eq:3ep1} \\
    (\eta \La_2 + \lambda \La_1 + \La_0 + \mu \Ma)\mathbf{v}                             & = \mathbf{0},\label{eq:3ep2} \\
    (\eta \mathbf{C}_2 + \lambda \mathbf{C}_1 + \mathbf{C}_0) \mathbf{w}                 & = \mathbf{0}\label{eq:3ep3}
  \end{align}
\end{subequations}
with
\begin{equation}
  \mathbf{C}_2 = \begin{bmatrix} 1 & 0 \\ 0& 0 \end{bmatrix},\
  \mathbf{C}_1 = \begin{bmatrix} 0 & 1 \\ 1 & 0 \end{bmatrix},\
  \mathbf{C}_0 = \begin{bmatrix} 0 & 0 \\ 0 & 1 \end{bmatrix}.
\end{equation}
It can be easily verified that \eqref{eq:3ep3} incorporates the equation
\begin{equation}
  \eta = \lambda^2,
\end{equation}
because a solution to this eigenvalue problem satisfies
\begin{equation}
  \det(\eta \mathbf{C}_2 + \lambda \mathbf{C}_1 + \mathbf{C}_0) = \eta - \lambda^2 = 0.
\end{equation}
The system \eqref{eq:3PEVP} is similar to those discussed in more detail in \cite{Kiefer2022, Gravenkamp2026}, and we use the MATLAB implementation available at \cite{Plestenjak2023} to solve it.
To briefly summarize the essential steps, we find solutions to the three-parameter EVP by solving related generalized eigenvalue problems (GEPs), see, e.g., \cite{Plestenjak2015} and the previous applications \cite{Kiefer2022, Gravenkamp2026}. These GEPs are
\begin{equation}\label{eq:delta_system}
  \mathbf{\Delta}_1 \mathbf{z} = \lambda \mathbf{\Delta}_0 \mathbf{z},\quad
  \mathbf{\Delta}_M \mathbf{z} = \mu \mathbf{\Delta}_0 \mathbf{z},\quad
  \mathbf{\Delta}_2 \mathbf{z} = \eta \mathbf{\Delta}_0 \mathbf{z},
\end{equation}
with the eigenvector $\mathbf{z}$ obtained by the Kronecker product $\mathbf{z}=\mathbf{u}\otimes \mathbf{v}\otimes \mathbf{w}$. The matrices in the above equations are referred to as operator determinants and are computed using the Kronecker product as:
\begin{alignat}{2}\label{eq:delta01}
  \mathbf{\Delta}_0 =   &\left|\begin{matrix}
                                  \mathbf{L}_2 & \mathbf{L}_1 & \mathbf{M} \cr
                                  \La_2        & \La_1        & \Ma \cr
                                  \mathbf{C}_2 & \mathbf{C}_1 & \mathbf{0}
                                \end{matrix}\right|_\otimes\negthickspace\negthickspace,\quad
  \mathbf{\Delta}_1 = -&\left|\begin{matrix}
                               \mathbf{L}_2 & \mathbf{L}_0 & \mathbf{M}\cr
                               \La_2        & \La_0        & \Ma \cr
                               \mathbf{C}_2 & \mathbf{C}_0 & \mathbf{0}
                             \end{matrix}\right|_\otimes
  \negthickspace\negthickspace,         \nonumber   \\
  \mathbf{\Delta}_2 = -&\left|\begin{matrix}
                                  \mathbf{L}_1 & \mathbf{L}_0 & \mathbf{M}\cr
                                  \La_1        & \La_0        & \Ma \cr
                                  \mathbf{C}_1 & \mathbf{C}_0 & \mathbf{0}
                                \end{matrix}\right|_\otimes
  \negthickspace\negthickspace,\quad
  \mathbf{\Delta}_M= -&\left|\begin{matrix}
                              \mathbf{L}_2 & \mathbf{L}_1 & \mathbf{L}_0 \cr
                              \La_2        & \La_1        & \La_0 \cr
                              \mathbf{C}_2 & \mathbf{C}_1 & \mathbf{C}_0
                            \end{matrix}\right|_\otimes \negthickspace\negthickspace.
\end{alignat}
Solving the first equation in \eqref{eq:delta_system} yields \textit{candidates} for $\eta$, $\lambda$, $\mu$ corresponding to critical points. Recovering the wavenumber $k = -\I\lambda$ and substituting into the original parameter-dependent eigenvalue problem \eqref{eq:evp_discrete} provides the corresponding frequencies $\omega$. It must be noted that \eqref{eq:delta_system} typically admits many additional solutions that are not critical points. Hence, we employ a simple postprocessing step, in which we select solutions with finite frequency that satisfy $c_p\approx c_g$ to a given tolerance. Specifically, if $k_\mathrm{cand}$, $\omega_\mathrm{cand}$ are candidate solutions, we consider them a valid critical point if $\omega_\mathrm{cand}$ is finite and 
\begin{equation}
  \left|\frac{|c_{p,\mathrm{cand}}|}{|c_{g,\mathrm{cand}}|}-1\right| < 10^{-4}
\end{equation}
with $c_{p,\mathrm{cand}} = \frac{\omega_\mathrm{cand}}{k_\mathrm{cand}}$ and $c_{g,\mathrm{cand}}$ computed using Eq.~\eqref{eq:groupVel}.

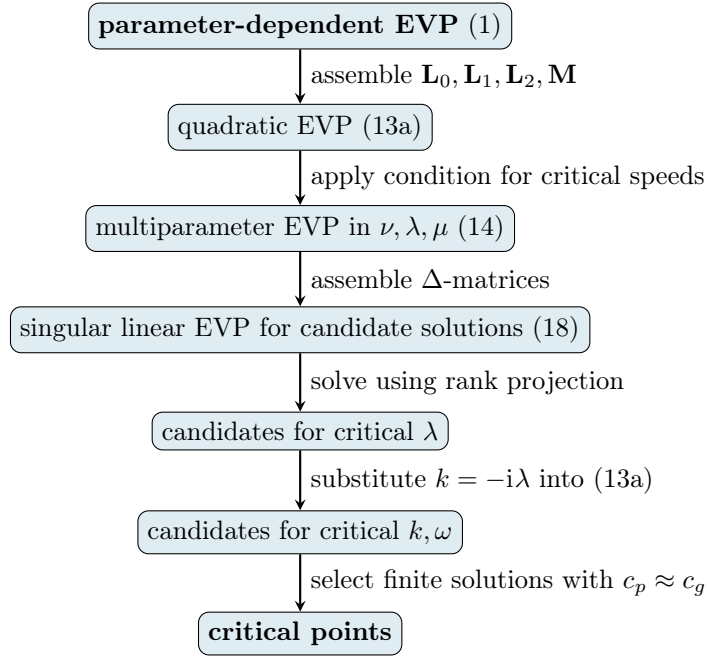
\begin{figure}
  \centering
\tikzstyle{terminator} = [rectangle, draw, text centered, rounded corners, minimum height=1em]
\tikzstyle{process} = [rectangle, draw, text centered, rounded corners, minimum height=1em]
\tikzstyle{decision} = [diamond, draw, text centered, minimum height=1em]
\tikzstyle{data}=[trapezium, draw, text centered, trapezium left angle=60, trapezium right angle=120, minimum height=1em]
\tikzstyle{connector} = [draw, -latex']
\begin{tikzpicture}[node distance = 2.8\baselineskip]
\node [terminator, fill=myblue!20] (gov) {\small \textbf{parameter-dependent EVP} \eqref{eq:evp_discrete}};
\node [process, fill=myblue!20, below of=gov] (QEP) {\small quadratic EVP \eqref{eq:2ep1}};
\node [process, fill=myblue!20, below of=QEP] (MEP) {\small multiparameter EVP in $\nu, \lambda,  \mu$ \eqref{eq:3PEVP}};
\node [process, fill=myblue!20, below of=MEP] (genEig) {\small singular linear EVP for candidate solutions \eqref{eq:delta_system}};
\node [process, fill=myblue!20, below of=genEig] (cand) {\small candidates for critical $\lambda$};
\node [process, fill=myblue!20, below of=cand] (candkw) {\small candidates for critical $k,\omega$};
\node [terminator, fill=myblue!20, below of=candkw] (sol) {\small \textbf{critical points}};
\draw [arrow] (gov) -- node[anchor=west] {\small assemble $\mathbf{L}_0$,\,$\mathbf{L}_1$,\,$\mathbf{L}_2$,\,$\mathbf{M}$ } (QEP);
\draw [arrow] (QEP) -- node[anchor=west] {\small apply condition for critical speeds} (MEP);
\draw [arrow] (MEP) -- node[anchor=west] {\small assemble $\Delta$-matrices} (genEig);
\draw [arrow] (genEig) -- node[anchor=west] {\small solve using rank projection} (cand);
\draw [arrow] (cand) -- node[anchor=west] {\small substitute $k=-\I\lambda$ into \eqref{eq:2ep1} } (candkw);
\draw [arrow] (candkw) -- node[anchor=west] {\small select finite solutions with $c_p\approx c_g$ } (sol);
\end{tikzpicture}
\caption{Summary of the essential steps in computing critical speeds. \label{fig:flowChart}}
\end{figure}

Since the GEPs are singular, the MATLAB toolbox MultiParEig \cite{Plestenjak2023} uses a rank-projection algorithm as proposed in \cite{Hochstenbach2023a}.
Due to their construction via Kronecker products, the size $n_\mathbf{\Delta}$ of the operator determinants increases quadratically with $n$ (the size of the original parameter-dependent eigenvalue problem). Specifically, it is given by the product of the sizes of the three coupled eigenvalue problems in \eqref{eq:3PEVP}, in our case  $n_\mathbf{\Delta} = n\negthinspace\cdot\negthinspace 2n\negthinspace\cdot\negthinspace 2 = 4n^2$. Computing all eigenvalues of the GEPs \eqref{eq:delta_system} is feasible on a current laptop computer for matrix sizes $n$ up to about 50, which already corresponds to $n_\mathbf{\Delta} = 10000$. In Fig.~\ref{fig:flowChart}, we present an overview of the essential steps in computing the critical points based on the proposed approach.

\subsection{Special case of a quartic eigenvalue problem}\noi
In principle, the approach outlined above can be straightforwardly extended to polynomial eigenvalue problems of arbitrary order (though the computational costs increase rapidly with the number of parameters). We present a special case of a fourth-order eigenvalue problem given by Eq.~\eqref{eq:operator_discrete_order4}, as it is encountered in the well-known formulation of a beam on an elastic foundation.
Differentiating the eigenvalue problem analogously and introducing $\xi = k^4$, $\mu = \omega^2$ leads to the two-parameter eigenvalue problem for critical points:
\begin{subequations}
  \label{eq:2PEVP_k4}
  \begin{align}
    (\xi \mathbf{L}_4 + \mathbf{L}_0 + \mu \mathbf{M})\mathbf{u} & = \mathbf{0},\label{eq:2ep1_k4} \\
    (\xi \La_4 + \La_0 + \mu \Ma)\mathbf{v}                      & = \mathbf{0}\label{eq:2ep2_k4}
  \end{align}
\end{subequations}
with
\begin{equation}
  \La_4 =
  \begin{bmatrix}
    \mathbf{L}_4  & \mathbf{0}   \\
    4\mathbf{L}_4 & \mathbf{L}_4
  \end{bmatrix}
\end{equation}
and $\La_0$, $\Ma$ defined as in \eqref{eq:criticalPointSystemExplicit}. Note that this simpler case leads to a two-parameter eigenvalue problem, as the original problem is already linear in $k^4$. Consequently, the operator determinants are of size $n_\mathbf{\Delta} = 2n^2$.

\section{Numerical experiments}\noi

\subsection{Illustrative example}\label{sec:acad}
\begin{figure}
  \subfloat[frequency]{\includegraphics[width = 0.49\textwidth]{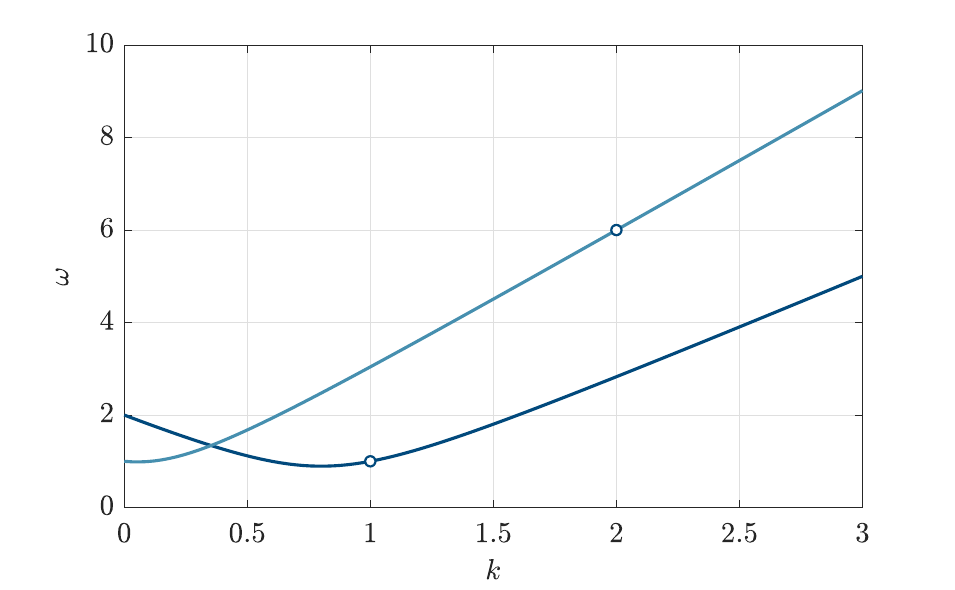}}\hfill
  \subfloat[phase and group velocity]{\includegraphics[width = 0.49\textwidth]{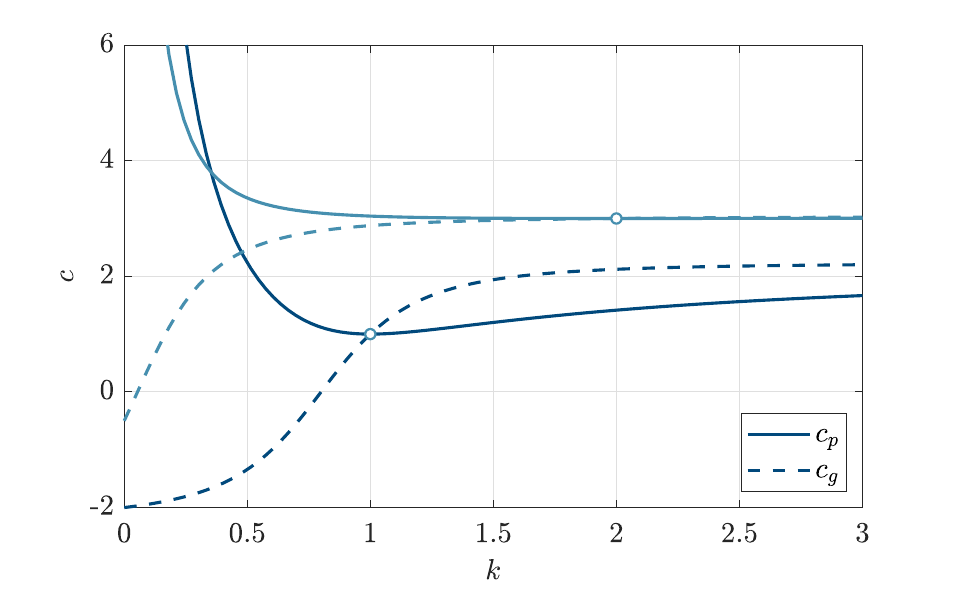}}
  \caption{Eigencurves of the minimal example in terms of frequency $\omega(k)$, phase velocity $c_p(k)$, and group velocity $c_g(k)$ (positive branches only). The critical points (indicated by `$\circ$') are the intersections of the phase velocity- and group velocity-curves. \label{fig:minimalExample}}
\end{figure}

\noi
We begin with an academic example featuring a simple solution that is easy to verify. Consider the parameter-dependent eigenvalue problem
\begin{equation}
  \label{eq:evp_minimalExample}
  (k^2 \mathbf{L}_2 + k \mathbf{L}_1 + \mathbf{L}_0 + \omega^2 \mathbf{M})\, \mathbf{u} = \mathbf{0}
\end{equation}
with
\begin{equation}
  \mathbf{L}_2 = -\begin{bmatrix} 153 & 64 \\ 64 & 57 \end{bmatrix},\quad
  \mathbf{L}_1 = \begin{bmatrix} 24 & -8 \\ -8 & 36 \end{bmatrix},\quad
  \mathbf{L}_0 = -\begin{bmatrix} 20 & 0 \\ 0 & 20 \end{bmatrix},\quad
  \mathbf{M} = \begin{bmatrix} 17 & 6 \\ 6 & 8 \end{bmatrix}.
\end{equation}
The matrix operator is of the form \eqref{eq:operator_discrete}, except that we define it directly in $k$ rather than $\I k$ for conciseness.
Solutions are given by the four eigencurves (including the positive and negative signs):
\begin{subequations}
  \label{eq:minimalEigencurves}
  \begin{align}
    \omega_1 & = \pm \sqrt{5\,k^2 - 8\,k + 4}, \\
    \omega_2 & = \pm \sqrt{9.25\,k^2 - k + 1}.
  \end{align}
\end{subequations}
The modes each exhibit one critical point at $k_1 = 1$ and $k_2 = 2$, respectively, i.e.,
\begin{equation}
  \omega_1'(k_1) = \frac{\omega_1(k_1)}{k_1} = 1,\quad \omega_2'(k_2) = \frac{\omega_2(k_2)}{k_2} = 3.
\end{equation}
Figure \ref{fig:minimalExample} shows the analytical eigencurves in terms of $\omega(k)$, $c_p(k)$, and $c_g(k)$, together with the critical points (listed in Tab.~\ref{tab:acad}), computed using the proposed approach. In this simple problem, the critical values are obtained exactly to machine precision.

\begin{table} \centering
  \caption{Critical points computed in example \ref{sec:acad}. \label{tab:acad}}
  \renewcommand{\arraystretch}{1} 
  \begin{tabular}{lrrr}\Xhline{2\arrayrulewidth}
      & $\omega_c $ & $k_c$  & $c_c$  \\
    \Xhline{2\arrayrulewidth}
    1 & 1.0    & 1.0 & 1.0 \\
    2 & 3.0    & 2.0 & 1.5 
    \\\Xhline{2\arrayrulewidth}
  \end{tabular}
\end{table}

\subsection{Beam on an elastic foundation}
\begin{figure}
  \subfloat[frequency]{\includegraphics[width = 0.49\textwidth]{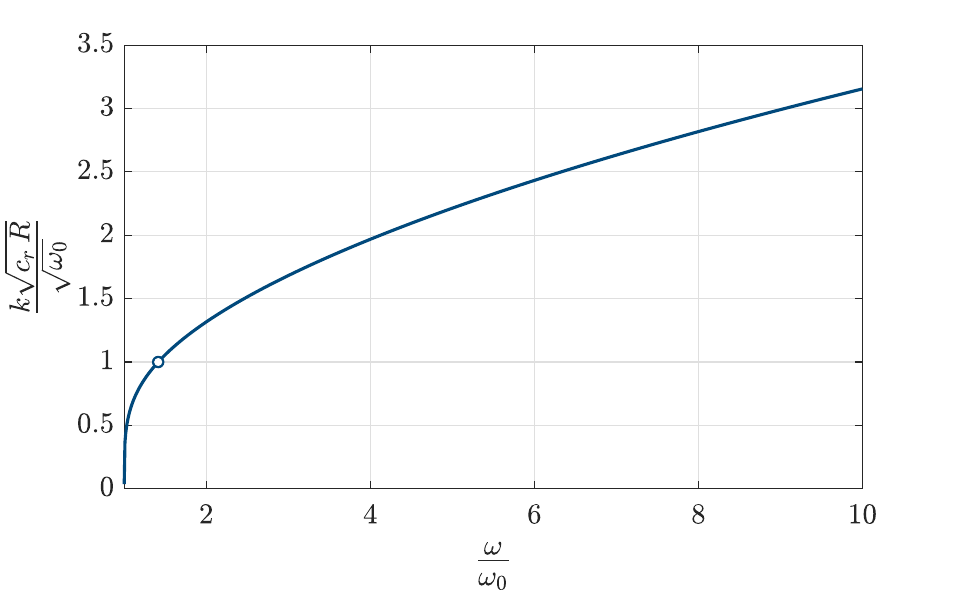}}\hfill
  \subfloat[phase and group velocity]{\includegraphics[width = 0.49\textwidth]{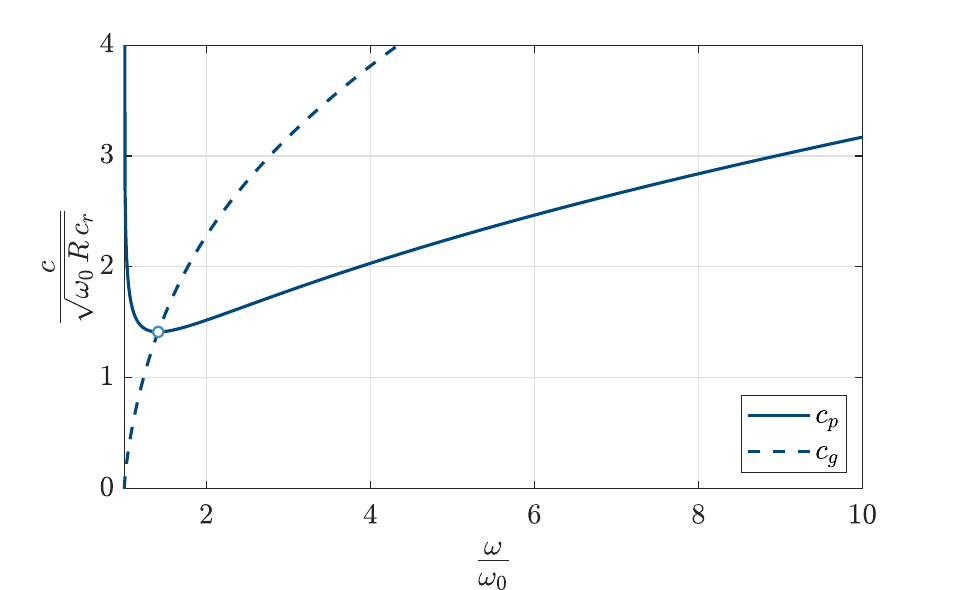}}
  \caption{Dimensionless eigencurves of the beam on an elastic foundation in terms of wavenumber $k(\omega)$, phase velocity $c_p(\omega)$, and group velocity $c_g(\omega)$ (positive branches only). The critical points ($\circ$) are the intersections of the phase velocity- and group velocity-curves. \label{fig:beam}}
\end{figure}

\noi
Let us consider a classical scalar benchmark example: a beam of infinite length on an elastic foundation. This classical problem is often used as a model for rails on soil layers, and the critical velocities of such simplified systems are well understood, see, e.g., \cite{Kausel2020} and the references therein.
Dispersion curves in terms of $k$ and $\omega$ are obtained as solutions to the scalar problem
\begin{equation}\label{eq:beamDispersion}
  \omega^2 - \omega_0^2 - k^4c_r^2R^2 = 0
\end{equation}
with the definitions
\begin{equation}
  \omega_0^2 \coloneq \frac{K}{\rho A},\qquad c_r^2 \coloneq \frac{E}{\rho},\qquad R^2 \coloneq \frac{I}{A},
\end{equation}
where $E$, $A$, $\rho$, $I$ denote the beam's properties (Young's modulus, cross-sectional area, mass density, moment of inertia), and $K$ is the foundation's stiffness per length. While this scalar problem clearly does not require numerical methods to obtain the critical speed, it is instructive for demonstrating the proposed method. In our notation established in Section~\ref{sec:cCrit}, we treat \eqref{eq:beamDispersion} as the scalar version of \eqref{eq:operator_discrete_order4} and pose the corresponding multiparameter eigenvalue problem of the form \eqref{eq:2PEVP_k4} with
\begin{equation}
  L_4 = -c_r^2R^2,\quad L_0 = -\omega_0^2, \quad M = 1,
\end{equation}
\begin{equation}
  \La_4 = -c_r^2R^2 \begin{bmatrix} 1 & 0 \\ 4 & 1 \end{bmatrix},\quad
  \La_0 = -\omega_0^2\begin{bmatrix} 1 & 0 \\ 0 & 1 \end{bmatrix},\quad
  \Ma = \begin{bmatrix} 1 & 0 \\ 2 & 1 \end{bmatrix}.
\end{equation}
The operator determinants are obtained as
\begin{equation}
  \mathbf{\Delta}_0 = -c_r^2R^2\begin{bmatrix} 0 & 0 \\ 2 & 0 \end{bmatrix},\quad
  \mathbf{\Delta}_1 = -\omega_0^2\begin{bmatrix} 0 & 0 \\ 2 & 0 \end{bmatrix},\quad
  \mathbf{\Delta}_2 = -c_r^2R^2\omega_0^2\begin{bmatrix} 0 & 0 \\ 4 & 0 \end{bmatrix}.
\end{equation}
In this simple case, it can be immediately observed that the operator determinants are, in fact, singular.
Wavenumbers and frequencies corresponding to critical points satisfy the eigenvalue problems
\begin{equation}
  \mathbf{\Delta}_1 \mathbf{x} = k_c^4 \mathbf{\Delta}_0 \mathbf{x},\quad \mathbf{\Delta}_2 \mathbf{x} = \omega_c^2 \mathbf{\Delta}_0 \mathbf{x}.
\end{equation}
The above eigenvalue problems have exactly one finite solution at 
\begin{equation}
k_c^4 = \frac{\omega_0^2}{c_r^2\,R^2},\qquad \omega_c^2=2\omega_0^2,\qquad  c_{c}=\frac{\omega_c}{k_c}
=\sqrt{2\,\omega_0\,c_r R}.
\end{equation}
The dispersion curves obtained by computing wavenumbers $k$ at discrete values of the frequency $\omega$, as well as the critical point obtained by solving the multiparameter eigenvalue problem by virtue of the proposed approach, are presented in Fig.~\ref{fig:beam}. The results are consistent with those presented in~\cite{Kausel2020}.

\subsection{Layered medium}\label{sec:soil}
\begin{figure}
  \subfloat[frequency]{\includegraphics[width = 0.48\textwidth]{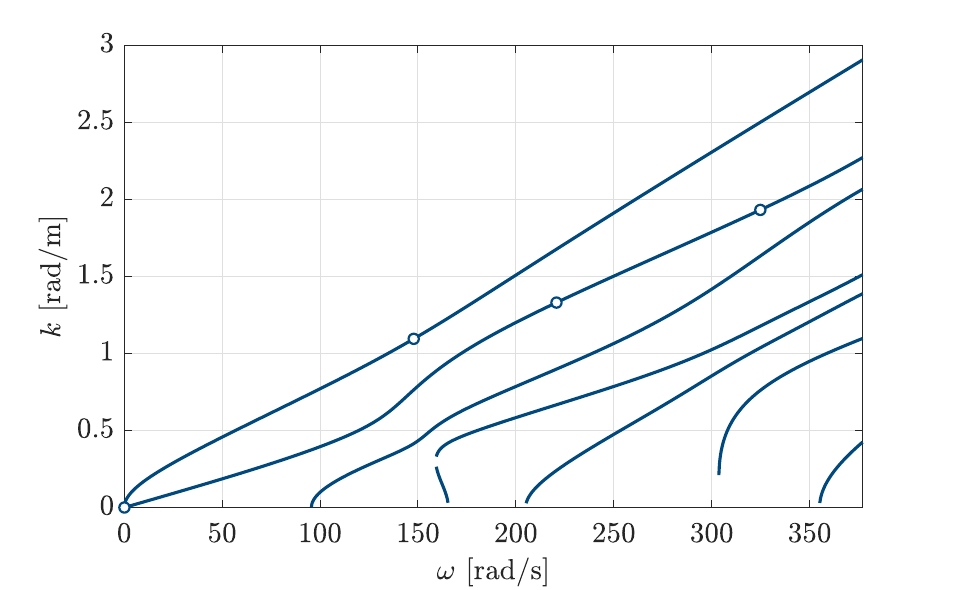}}\hfill
  \subfloat[phase and group velocity]{\includegraphics[width = 0.48\textwidth]{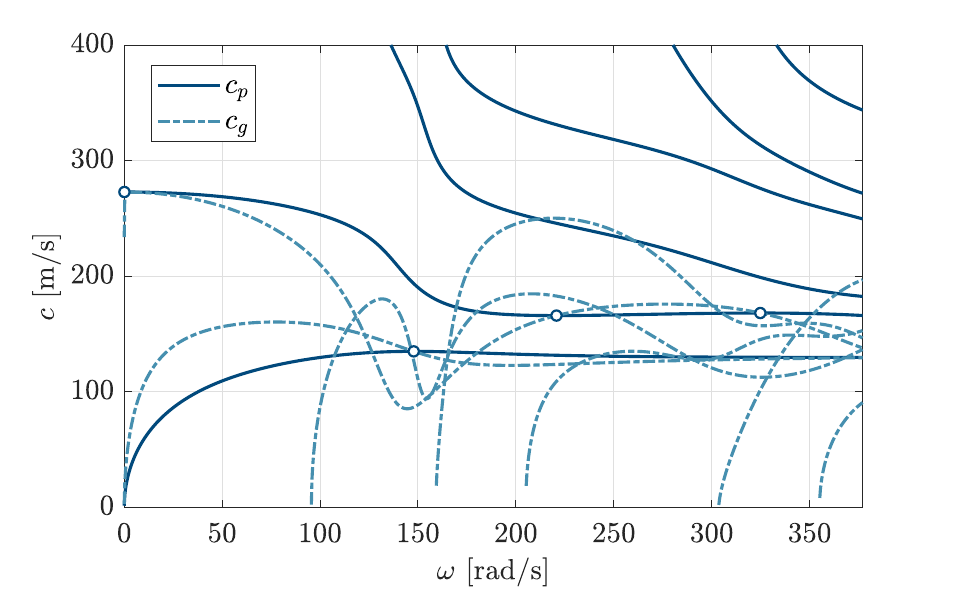}}
  \caption{Eigencurves of the layered soil in terms of wavenumber $k(\omega)$, phase velocity $c_p(\omega)$, and group velocity $c_g(\omega)$. The critical points ($\circ$) are the intersections of the phase velocity- and group velocity-curves. \label{fig:soil}}
\end{figure}

\noi
Moving to a more realistic, practical case, this example is adapted from that presented by Kausel \cite{Kausel2020}, who studied critical velocities in stratified infinite media representing the foundation beneath rail systems. Specifically, we consider a system consisting of two layers of thicknesses 2\,m and 3\,m, respectively, representing ballast and embankment with the material parameters listed in Table~\ref{tab:materialParameters} (in terms of shear wave velocity $c_s$, mass density $\rho$, and Poisson's ratio $\nu$). The two layers are fully coupled (i.e., assuming continuous displacements across the interface), and the outer boundaries are traction-free.
Following the semi-analytical approach described in \cite{Gravenkamp2012e}, each of the two layers is discretized by only one higher-order finite element. For the presented frequency range up to $60\,\unit{Hz}$ ($377\,\unitfrac{rad}{s}$), an element order of five is found to be sufficient, leading to a total of 22 degrees of freedom in a plane strain model. The matrices and dispersion curves have been computed using the open-source MATLAB code \emph{SAMWISE} \cite{Gravenkamp2024}. Results are presented in Fig.~\ref{fig:soil}, showing three critical points. Interestingly, two of the critical points belong to the same mode and lie on a branch with small curvature. The numerical values of the three critical points are listed in Table~\ref{tab:soil_critical}.

\begin{table} \centering
  \caption{Material parameters used in example \ref{sec:soil}: shear wave velocity $c_s$, mass density $\rho$, and Poisson's ratio $\nu$, thickness $d$. \label{tab:materialParameters}}
  \renewcommand{\arraystretch}{1} 
  \begin{tabular}{clcccc}\Xhline{2\arrayrulewidth}
    layer & material   & $c_s$\,[$\unitfrac{m}{s}$] & $\rho$\,[$\unitfrac{kg}{m^3}$] & $\nu$ & $d$ [m] \\
    \Xhline{2\arrayrulewidth}
    1     & ballast    & 200                        & 2000                           & 0.25  & 2       \\
    2     & embankment & 141                        & 2000                           & 0.25  & 3
    \\\Xhline{2\arrayrulewidth}
  \end{tabular}
\end{table}

\begin{table} \centering
  \caption{Critical points computed in example \ref{sec:soil}. \label{tab:soil_critical}}
  \renewcommand{\arraystretch}{1} 
  \begin{tabular}{lrrr}\Xhline{2\arrayrulewidth}
      & $\omega_c\,[\unitfrac{rad}{s}]$ & $k_c\,[\unitfrac{rad}{m}]$ & $c_c\,[\unitfrac{m}{s}]$  \\
    \Xhline{2\arrayrulewidth}
    1 & 147.83    & 1.09 & 135.08 \\
    2 & 220.73    & 1.33 & 166.02 \\
    3 & 324.82    & 1.93 & 168.24
    \\\Xhline{2\arrayrulewidth}
  \end{tabular}
\end{table}

\subsection{Three-dimensional structure}\label{sec:ex3d}
\begin{figure}\centering
  \includegraphics[width = 0.6\textwidth]{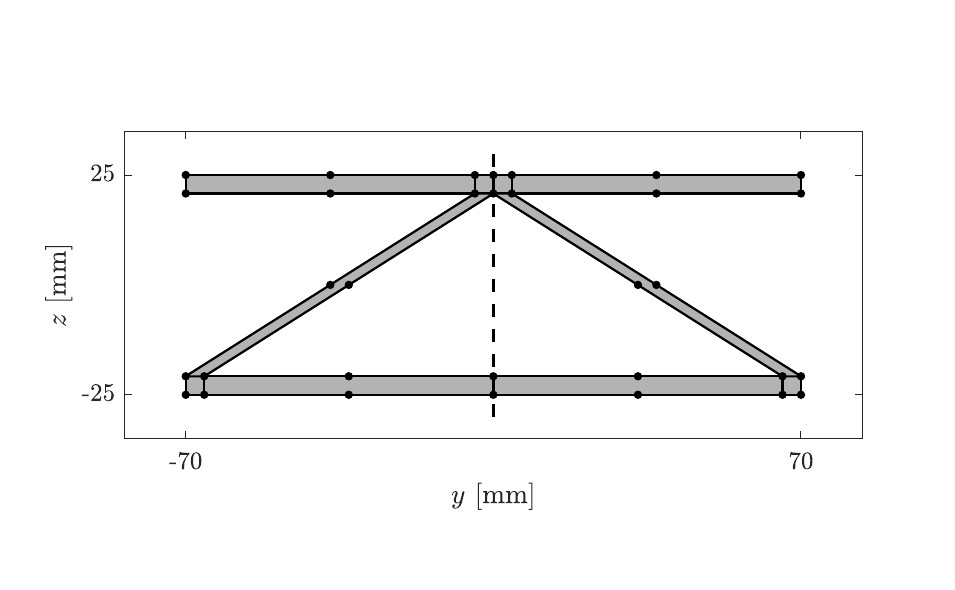}
  \caption{Cross-section geometry of the structure studied in Section~\ref{sec:ex3d}. \label{fig:panel_geometry}}
\end{figure}

\begin{figure}
  \subfloat[frequency]{\includegraphics[width = 0.48\textwidth]{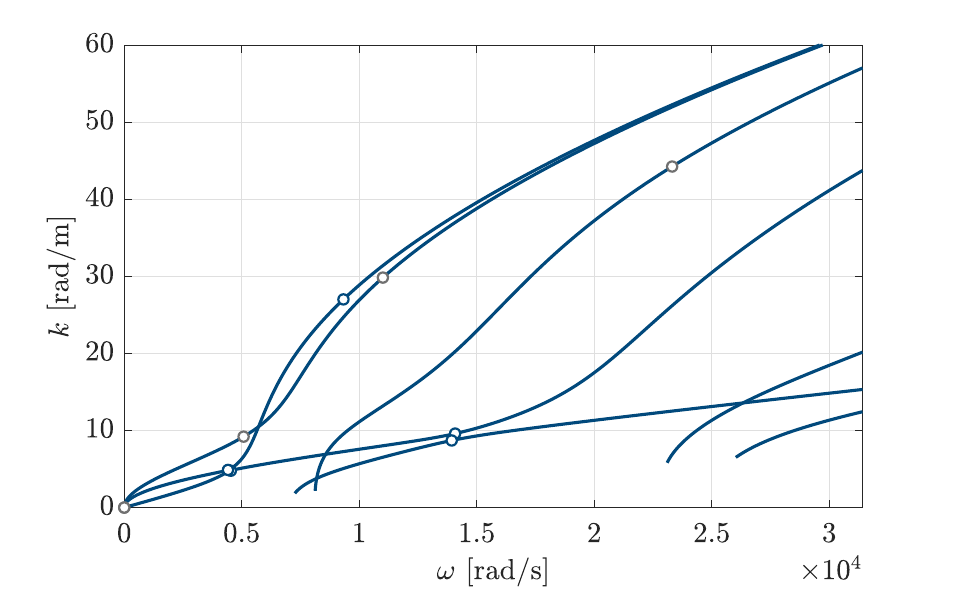}}\hfill
  \subfloat[phase and group velocity]{\includegraphics[width = 0.48\textwidth]{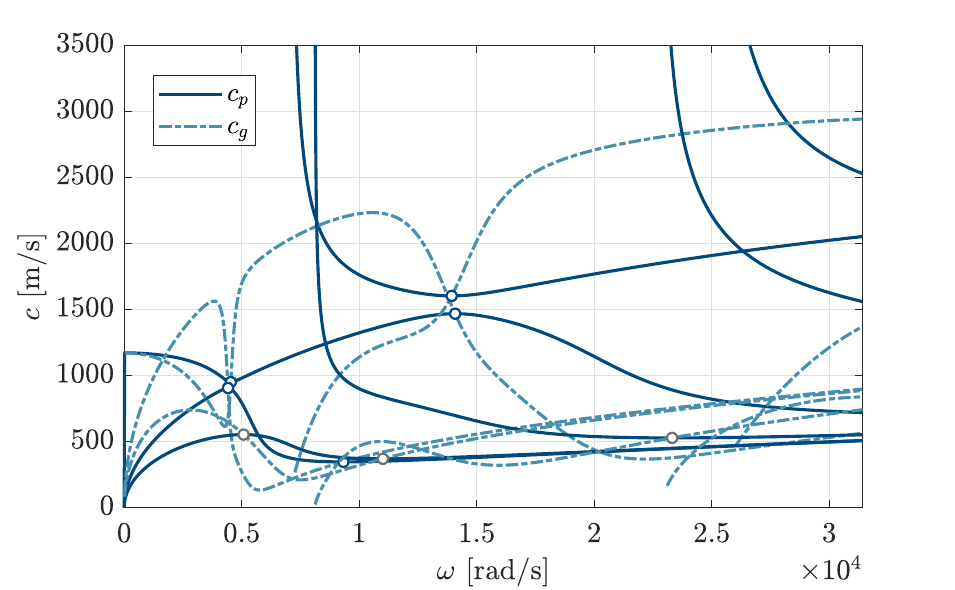}}
  \caption{Eigencurves of the structure depicted in Fig.~\ref{fig:panel_geometry} in terms of wavenumber $k(\omega)$, phase velocity $c_p(\omega)$, and group velocity $c_g(\omega)$. The critical points ($\circ$) are the intersections of the phase velocity- and group velocity-curves. \label{fig:panelCurves}}
\end{figure}

\noi 
As a final test, we consider a thin-walled prismatic waveguide whose cross-section is depicted in Fig.~\ref{fig:panel_geometry}. The structure is assumed to be of infinite extent in the $x$-direction. Thin-walled built-up members of this type are representative of lightweight transport structures, where both in-plane and out-of-plane deformation mechanisms coexist. This example is therefore more demanding than the preceding layered-medium benchmark: it involves symmetries, dispersion branches interact closely, and multiple stationary points of the phase velocity occur within a moderate frequency range. Hence, the structure provides a useful test case for assessing whether the proposed multiparameter formulation can identify critical points without relying on branch tracing. The total width and height of the cross-section are 140\,mm and 50\,mm, respectively, and the thickness of each element is 4.2\,mm. The material is defined by a Young's modulus of 71\,GPa, a mass density of 2700\,kg/m$^3$, and a Poisson's ratio of 0.332.
As the current approach requires the matrix size of the discretized model to be small, we use a coarse discretization, which nevertheless gives sufficiently accurate results for frequencies up to 5\,kHz ($10^4\,\uppi\,\unitfrac{rad}{s}$). Since the cross-section is symmetric with respect to the $y$-axis, we discretize half of it and apply, in turn, symmetric and antisymmetric boundary conditions (see \cite{Gravenkamp2013e} for implementation details). Half of the cross-section is discretized by five finite elements. Note that the elements use a 'mixed-order-approach' with linear interpolation across the thickness direction and quadratic interpolation along each member of the cross-section. This discretization leads to 18 nodes and, consequently, 54 degrees of freedom (50 and 46 after applying symmetric and antisymmetric boundary conditions, respectively). Again, the matrices and dispersion curves were computed using the software \emph{SAMWISE} \cite{Gravenkamp2024}. Results are presented in Fig.~\ref{fig:panelCurves}. This structure exhibits eight critical points within the selected frequency range. Numerical values of the computed critical points are listed in Table~\ref{tab:3d_critical}.

\begin{table} \centering
  \caption{Critical points computed in example \ref{sec:ex3d}. \label{tab:3d_critical}}
  \renewcommand{\arraystretch}{1} 
  \begin{tabular}{lrrr}\Xhline{2\arrayrulewidth}
      & $\omega\,[\unitfrac{rad}{s}]$ & $k\,[\unitfrac{rad}{m}]$ & $c_p\,[\unitfrac{m}{s}]$ \\
    \Xhline{2\arrayrulewidth}
    1 & 4416.4                        & 4.88                     & 904.59                   \\
    2 & 4536.1                        & 4.78                     & 949.46                   \\
    3 & 5080.1                        & 9.20                     & 552.27                   \\
    4 & 9329.1                        & 27.00                    & 345.47                   \\
    5 & 11007.0                       & 29.83                    & 368.95                   \\
    6 & 13935.9                       & 8.70                     & 1601.78                  \\
    7 & 14080.9                       & 9.60                     & 1467.13                  \\
    8 & 23315.8                       & 44.24                    & 527.02
    \\\Xhline{2\arrayrulewidth}
  \end{tabular}
\end{table}

\section{Conclusion}\noi
We have presented a direct approach to computing critical velocities in elastic waveguides. By incorporating the condition of coinciding phase and group velocities into the parameter-dependent eigenvalue problem arising in typical semi-analytical approaches, the identification of critical points is reformulated as a multiparameter eigenvalue problem. This formulation avoids the need for dispersion-curve tracing or complicated postprocessing.

We employed previously established approaches and algorithms to solve the resulting singular polynomial multiparameter eigenvalue problems in a straightforward manner. The numerical examples demonstrate that the proposed approach reliably identifies critical points in different scenarios, including simple benchmark systems, layered media, and three-dimensional prismatic structures. In particular, the method remains effective in situations with multiple critical points, mode crossings, and nearly flat (low-curvature) minima in the dispersion curves. Currently, the main drawback of the proposed approach is that its computational cost increases rapidly with matrix size, making it de facto feasible only for small systems. Hence, future work will focus on improving computational efficiency by assessing strategies similar to those proposed for computing ZGV points \cite{Kiefer2022, Plestenjak2025}, as well as on extending the framework to more general settings, including periodically structured media and leaky-wave problems.

\bibliographystyle{elsarticle-num}
\bibliography{criticalVelocity.bib}

@article{Bierer2007,
  title = {A semi-analytical model in time domain for moving loads},
  author = {Bierer, T and Bode, C},
  year = 2007,
  month = dec,
  journal = {Soil Dynamics and Earthquake Engineering},
  volume = {27},
  number = {12},
  pages = {1073--1081},
  issn = {02677261},
  doi = {10.1016/j.soildyn.2007.03.008}
}

@article{Estaire2024,
  title = {On the use of dispersion curves to determine the critical speed of railway tracks. Application to case studies},
  author = {Estaire, Jos{\'e} and {Crespo-Chac{\'o}n}, In{\'e}s},
  year = 2024,
  journal = {Transportation Geotechnics},
  volume = {46},
  pages = {101225},
  issn = {22143912},
  doi = {10.1016/j.trgeo.2024.101225}
}

@article{Faragau2025,
  title = {Controlling instability of high-speed magnetically suspended vehicles: The interaction of the electromagnetic and wave-induced instability mechanisms},
  shorttitle = {Controlling instability of high-speed magnetically suspended vehicles},
  author = {F{\u a}r{\u a}g{\u a}u, Andrei B. and Metrikine, Andrei V. and Paul, Jithu and Van Leijden, Rens and Van Dalen, Karel N.},
  year = 2025,
  journal = {Journal of Sound and Vibration},
  volume = {608},
  pages = {119077},
  issn = {0022460X},
  doi = {10.1016/j.jsv.2025.119077}
}

@article{Gravenkamp2012e,
  title = {A numerical approach for the computation of dispersion relations for plate structures using the scaled boundary finite element method},
  author = {Gravenkamp, Hauke and Song, Chongmin and Prager, Jens},
  year = 2012,
  journal = {Journal of Sound and Vibration},
  volume = {331},
  number = {11},
  pages = {2543--2557},
  publisher = {Elsevier},
  doi = {10.1016/j.jsv.2012.01.029}
}

@article{Gravenkamp2013e,
  title = {The computation of dispersion relations for three-dimensional elastic waveguides using the Scaled Boundary Finite Element Method},
  author = {Gravenkamp, Hauke and Man, Hou and Song, Chongmin and Prager, Jens},
  year = 2013,
  journal = {Journal of Sound and Vibration},
  volume = {332},
  pages = {3756--3771},
  publisher = {Elsevier},
  issn = {0022460X},
  doi = {10.1016/j.jsv.2013.02.007}
}

@article{Gravenkamp2018,
  title = {Efficient simulation of elastic guided waves interacting with notches, adhesive joints, delaminations and inclined edges in plate structures},
  author = {Gravenkamp, Hauke},
  year = 2018,
  journal = {Ultrasonics},
  volume = {82},
  pages = {101--113},
  publisher = {Elsevier B.V.},
  issn = {0041624X},
  doi = {10.1016/j.ultras.2017.07.019}
}

@article{Gravenkamp2021b,
  title = {High-order shape functions in the scaled boundary finite element method revisited},
  author = {Gravenkamp, Hauke and Saputra, Albert and Duczek, Sascha},
  year = 2021,
  journal = {Archives of Computational Methods in Engineering},
  volume = {28},
  pages = {473--494},
  issn = {1886-1784},
  doi = {10.1007/s11831-019-09385-1}
}

@article{Gravenkamp2023i,
  title = {Notes on osculations and mode tracing in semi-analytical waveguide modeling},
  author = {Gravenkamp, Hauke and Plestenjak, Bor and Kiefer, Daniel A.},
  year = 2023,
  journal = {Ultrasonics},
  volume = {135},
  pages = {107112},
  doi = {10.1016/j.ultras.2023.107112}
}

@misc{Gravenkamp2024,
  title = {SAMWISE: Semi-Analytical Modeling of Waves in Structural Elements},
  author = {Gravenkamp, Hauke},
  year = 2024,
  doi = {10.5281/zenodo.13830671}
}

@misc{Gravenkamp2024a,
  title = {Leaky Guided Waves Examples},
  author = {Gravenkamp, Hauke and Plestenjak, Bor and Kiefer, Daniel A. and Jarlebring, Elias},
  year = 2024,
  doi = {10.5281/zenodo.13825263}
}

@article{Gravenkamp2025,
  title = {Computation of leaky waves in layered structures coupled to unbounded media by exploiting multiparameter eigenvalue problems},
  author = {Gravenkamp, Hauke and Plestenjak, Bor and Kiefer, Daniel A and Jarlebring, Elias},
  year = 2025,
  journal = {Journal of Sound and Vibration},
  volume = {596},
  pages = {118716},
  doi = {10.1016/j.jsv.2024.118716}
}

@article{Gravenkamp2026,
  title = {On the direct numerical computation of Hopf bifurcations to assess the dynamic stability of fluid-conveying cantilevered pipes},
  author = {Gravenkamp, Hauke and Plestenjak, Bor},
  year = 2026,
  journal = {Computers \& Structures},
  volume = {320},
  pages = {108039},
  doi = {10.1016/j.compstruc.2025.108039}
}

@article{Hayashi2003a,
  title = {Guided wave dispersion curves for a bar with an arbitrary cross-section, a rod and rail example},
  author = {Hayashi, Takahiro and Song, W.-J. and Rose, J L},
  year = 2003,
  month = may,
  journal = {Ultrasonics},
  volume = {41},
  number = {3},
  pages = {175--183},
  issn = {0041624X},
  doi = {10.1016/S0041-624X(03)00097-0}
}

@article{Hochstenbach2023a,
  title = {Solving singular generalized eigenvalue problems. Part II: projection and augmentation},
  shorttitle = {Solving singular generalized eigenvalue problems. Part II},
  author = {Hochstenbach, Michiel E. and Mehl, Christian and Plestenjak, Bor},
  year = 2023,
  journal = {SIAM Journal on Matrix Analysis and Applications},
  volume = {44},
  number = {4},
  pages = {1589--1618},
  issn = {0895-4798, 1095-7162},
  doi = {10.1137/22M1513174}
}

@article{Kausel1981a,
  title = {Stiffness matrices for layered soils},
  author = {Kausel, Eduardo and Ro{\"e}sset, Jos{\'e} Manuel},
  year = 1981,
  journal = {Bulletin of the Seismological Society of America},
  volume = {71},
  number = {6},
  pages = {1743--1761},
  publisher = {Seismological Society of America},
  doi = {10.1785/BSSA0710061743}
}

@article{Kausel1994,
  title = {Thin-layer method: Formulation in the time domain},
  author = {Kausel, Eduardo},
  year = 1994,
  journal = {International Journal for Numerical Methods in Engineering},
  volume = {37},
  pages = {927--941},
  doi = {10.1002/nme.1620370604}
}

@article{Kausel2015,
  title = {Osculations of spectral lines in a layered medium},
  author = {Kausel, Eduardo and Malischewsky, Peter and Barbosa, Jo{\~a}o},
  year = 2015,
  journal = {Wave Motion},
  volume = {56},
  pages = {22--42},
  publisher = {Elsevier B.V.},
  issn = {01652125},
  doi = {10.1016/j.wavemoti.2015.01.004}
}

@article{Kausel2020,
  title = {Proof of critical speed of high-speed rail underlain by stratified media},
  author = {Kausel, Eduardo and Estaire, Jos{\'e} and {Crespo-Chac{\'o}n}, In{\'e}s},
  year = 2020,
  journal = {Proceedings of the Royal Society A},
  volume = {476},
  pages = {20200083}
}

@article{Kiefer2022,
  title = {Computing zero-group-velocity points in anisotropic elastic waveguides: globally and locally convergent methods},
  author = {Kiefer, Daniel A. and Plestenjak, Bor and Gravenkamp, Hauke and Prada, Claire},
  year = 2023,
  journal = {The Journal of the Acoustical Society of America},
  volume = {153},
  number = {2},
  pages = {1386--1398},
  doi = {10.1121/10.0017252}
}

@article{Krylov1995,
  title = {Generation of ground vibrations by superfast trains},
  author = {Krylov, Victor V.},
  year = 1995,
  journal = {Applied Acoustics},
  volume = {44},
  pages = {149--164},
  issn = {0003682X},
  doi = {10.1016/0003-682X(95)91370-I},
  copyright = {https://www.elsevier.com/tdm/userlicense/1.0/}
}

@article{Lamb1917,
  title = {On waves in an elastic plate},
  author = {Lamb, Horace},
  year = 1917,
  journal = {Proceedings of the Royal Society of London},
  volume = {93},
  number = {648},
  pages = {114--128},
  publisher = {The Royal Society},
  doi = {10.1098/rspa.1917.0008}
}

@article{Liu2022a,
  title = {Exact wave propagation analysis of lattice structures based on the dynamic stiffness method and the Wittrick--Williams algorithm},
  author = {Liu, Xiang and Lu, Zhaoming and Adhikari, Sondipon and Li, YingLi and Banerjee, J. Ranjan},
  year = 2022,
  journal = {Mechanical Systems and Signal Processing},
  volume = {174},
  pages = {109044},
  issn = {08883270},
  doi = {10.1016/j.ymssp.2022.109044}
}

@article{Liu2024b,
  title = {A wavenumber dynamic stiffness method for exact and efficient dispersion analysis of plate built-up waveguides},
  author = {Liu, Xiang and Zhou, Weixian and Filippi, Matteo and Wang, Yu},
  year = 2024,
  journal = {Journal of Sound and Vibration},
  volume = {591},
  pages = {118605},
  issn = {0022460X},
  doi = {10.1016/j.jsv.2024.118605}
}

@article{Liu2026,
  title = {Exact and efficient modal and critical speed analyses of high-speed railway catenary systems based on the dynamic stiffness method},
  author = {Liu, Xiang and Wang, Zekai and Sun, Chengli and Fang, Congcong and Lu, Tao and Zhao, Haibo},
  year = 2026,
  journal = {Computers \& Structures},
  volume = {327},
  pages = {108249},
  issn = {00457949},
  doi = {10.1016/j.compstruc.2026.108249}
}

@article{Marzani2008,
  title = {A semi-analytical finite element formulation for modeling stress wave propagation in axisymmetric damped waveguides},
  author = {Marzani, Alessandro and Viola, Erasmo and Bartoli, Ivan and {Lanza di Scalea}, Francesco and Rizzo, Piervincenzo},
  year = 2008,
  journal = {Journal of Sound and Vibration},
  volume = {318},
  pages = {488--505},
  issn = {0022460X},
  doi = {10.1016/j.jsv.2008.04.028}
}

@article{Mezher2016,
  title = {Railway critical velocity -- Analytical prediction and analysis},
  author = {Mezher, Sara B. and Connolly, David P. and Woodward, Peter K. and Laghrouche, Omar and Pombo, Joao and Costa, Pedro Alves},
  year = 2016,
  journal = {Transportation Geotechnics},
  volume = {6},
  pages = {84--96},
  issn = {22143912},
  doi = {10.1016/j.trgeo.2015.09.002}
}

@article{Muhic2010,
  title = {On the quadratic two-parameter eigenvalue problem and its linearization},
  author = {Muhi{\v c}, Andrej and Plestenjak, Bor},
  year = 2010,
  journal = {Linear Algebra and Its Applications},
  volume = {432},
  number = {10},
  pages = {2529--2542},
  issn = {0024-3795},
  doi = {10.1016/j.laa.2009.12.022}
}

@article{Plestenjak2015,
  title = {Spectral collocation for multiparameter eigenvalue problems arising from separable boundary value problems},
  author = {Plestenjak, Bor and Gheorghiu, C{\u a}lin I. and Hochstenbach, Michiel E.},
  year = 2015,
  journal = {Journal of Computational Physics},
  volume = {298},
  pages = {585--601},
  issn = {10902716},
  doi = {10.1016/j.jcp.2015.06.015}
}

@misc{Plestenjak2023,
  title = {MultiParEig},
  author = {Plestenjak, Bor},
  year = 2023,
  url = {www.mathworks.com/matlabcentral/fileexchange/47844-multipareig}
}

@article{Plestenjak2025,
  title = {A Sylvester equation approach for the computation of zero-group-velocity points in waveguides},
  author = {Plestenjak, Bor and Kiefer, Daniel A. and Gravenkamp, Hauke},
  year = 2025,
  journal = {Computational Mechanics},
  doi = {10.1007/s00466-025-02656-8}
}

@article{Tran2025,
  title = {Analysis of dynamic response of infinite beam on a periodical viscoelastic foundation subjected to moving loads and calculation of the critical train speed},
  author = {Tran, Le-Hung and Le, Thuy-Duong and Schmidt, Franziska},
  year = 2025,
  journal = {Archive of Applied Mechanics},
  volume = {95},
  number = {172},
  pages = {1--17},
  issn = {0939-1533, 1432-0681},
  doi = {10.1007/s00419-025-02871-y}
}

@article{Vesali2021,
  title = {Analysis of conceptual similarities and differences of wave speed and critical speed in the overhead catenary system},
  author = {Vesali, Farzad and Rezvani, Mohammad Ali and Molatefi, Habibolah},
  year = 2021,
  journal = {Measurement},
  volume = {176},
  pages = {109164},
  issn = {02632241},
  doi = {10.1016/j.measurement.2021.109164}
}

@article{Xin2017,
  title = {Self-controlled wave propagation in hyperelastic media},
  author = {Xin, Fengxian and Lu, Tian Jian},
  year = 2017,
  journal = {Scientific Reports},
  volume = {7},
  pages = {7581},
  issn = {2045-2322},
  doi = {10.1038/s41598-017-08098-4}
}

@article{Zhou2024,
  title = {Wavenumber dynamic stiffness formulation for exact dispersion analysis of moderately thick symmetric cross-ply laminated plate built-up waveguides},
  author = {Zhou, Weixian and Liu, Xiang and Wang, Yu and Zhao, Xueyi},
  year = 2024,
  journal = {Thin-Walled Structures},
  volume = {204},
  pages = {112305},
  issn = {02638231},
  doi = {10.1016/j.tws.2024.112305}
}

\end{document}